\documentclass[12pt,reqno,a4paper]{amsart}
\usepackage[colorlinks,linkcolor=blue]{hyperref}

\newcommand{\Span}[1]{\left<#1\right>}
\newcommand{\ga}{\gamma}
\newcommand{\la}{\lambda}
\newcommand{\om}{\omega}
\newcommand{\aff}{\mathbb A}
\newcommand{\PP}{\mathbb P}
\newcommand{\Oh}{\mathcal O}
\newcommand{\sI}{\mathcal I}
\newcommand{\iso}{\cong}
\newcommand{\into}{\hookrightarrow}
\newcommand{\bij}{\leftrightarrow}
\newcommand{\bijj}{\longleftrightarrow}
\DeclareMathOperator{\Der}{Der}
\DeclareMathOperator{\Grass}{Grass}
\DeclareMathOperator{\Hom}{Hom}
\DeclareMathOperator{\Pic}{Pic}
\DeclareMathOperator{\Segre}{Segre}

\title[Tutorial on Tom and Jerry]{Tutorial on Tom and Jerry:
\\ the two smoothings of the
\\ anticanonical cone over $\PP(1,2,3)$}

\author{Gavin Brown}

\author{Miles Reid}
\address{Mathematics Institute,
University of Warwick.
Coventry CV4 7AL, England}
\email{G.Brown@warwick.ac.uk}
\email{Miles.Reid@warwick.ac.uk}

\author{Jan Stevens}
\address{Department of Mathematical Sciences, Chalmers University of
Technology and University of Gothenburg.
SE 412 96 Gothenburg, Sweden}
\email{stevens@chalmers.se}

\date{}

\begin{document}
\begin{abstract}
This is a first introduction to unprojection methods, and more
specifically to Tom and Jerry unprojections. These two harmless tricks
deserve to be better known, since they answer many practical questions
about constructing codimension~4 Gorenstein subschemes. In particular,
we discuss here the two smoothing components of the anti\-canonical cone
over $\PP(1,2,3)$.
\end{abstract}

\maketitle

Section~\ref{s!1} treats the ``$6\times6$ extrasymmetric format'', that
describes the Segre embedding of $\PP^2\times\PP^2$ and some of its
degenerations. One can view this as just algebraic manipulations, or as
a typical case of Tom unprojection. In a similar vein, Section~\ref{s!2}
treats the ``Double Jerry construction'', that describes the Segre
embedding of $\PP^1\times\PP^1\times\PP^1$ and some of its
degenerations. In Section~\ref{s!3} we put these two unprojection
constructions together as a versal deformation of the anticanonical cone
over $\PP(1,2,3)$ over a reducible base, with the obstructions also
controlled by the matrix format. We conclude with some general remarks,
mnemonics, slogans, and FAQ. We do not pretend any generality, or any
theoretical treatment of Gorenstein codimension~4 (compare \cite{G4}).

\section{The anticanonical cone over $\PP(1,2,3)$}
\label{s!0}
Let $X\subset\aff^7$ be the anticanonical cone over
$\PP(1,2,3)_{\Span{u,v,w}}$; this is also the quotient by the group
action $\frac16(1,2,3)$ on $\aff^3_{\Span{u,v,w}}$. We set out its 7
coordinate monomials as the Newton polygon
\begin{equation}
\renewcommand{\arraystretch}{1.5}
\renewcommand{\arraycolsep}{0.4em}
\begin{matrix}
u^6 & u^4v & u^2v^2 & v^3 \\
u^3w & uvw \\
w^2
\end{matrix}
\qquad = \qquad
\renewcommand{\arraycolsep}{0.9em}
\begin{matrix}
a & b & c & x \\
d & e \\
f
\end{matrix}
\label{eq!1}
\end{equation}
The somewhat idiosyncratic choice of coordinates on $\aff^7$ relates
to the extrasymmetric format of Section~2.

One finds the equations defining $X$ without difficulty. The semigroup
ideal of internal monomials of the Newton polygon is generated by the
single monomial $e=uvw$. There are tag relations between any three
consecutive boundary monomials, that involve $e$ if we turn a corner:
\begin{equation}
ac-b^2,\quad
xb-c^2,\quad
cf-e^2,\quad
xdf-e^3,\quad
af-d^2,\quad
bd-ae. \notag
\end{equation}
Note in particular the equations $cf=e^2$ (that is, the tag at $x$ is
$0$) and $xd=f^{-1}e^3$ or $xdf=e^3$ (the tag at $f$ is $-1$).

These equations define the toric variety $X$ in the complement of the
coordinate hyperplanes, where $e$ is invertible. The remaining
generators of $I_X$ come by coloning out $e$: for example, $cf-e^2$ and
$xdf-e^3$ give $\big(c(xdf-e^3)-xd(cf-e^2)\big)e^{-2}=xd-ce$ where $e$ is
invertible. The ideal is generated by the 9 binomials:
\begin{equation} \label{eq!2}
\begin{gathered}
ac-b^2,\quad
xb-c^2,\quad
cf-e^2,\quad
af-d^2,\quad
bd-ae \\
xd-ce,\quad
bf-de,\quad
dc-be,\quad
xa-bc.
\end{gathered}
\end{equation}

Another way to view the equations is that they describe a singular del
Pezzo surface $S$ of degree 6. The monomials $U=u^3$, $V=uv$, $W=w$ base
$H^0(\PP(1,2,3),\Oh(3))$. We view them as coordinates on $\PP^2$. Then
multiplying \eqref{eq!1} by $u^3$ gives the 7 monomials
\begin{equation}
\renewcommand{\arraystretch}{1.8}
\renewcommand{\arraycolsep}{0.6em}
\begin{matrix}
U^3 & U^2V & UV^2 & V^3 \\
U^2W & UVW \\
UW^2
\end{matrix}
\label{eq!dP}
\end{equation}
that base the linear system of cubics in $\PP^2_{U,V,W}$ with flex line $U=0$ at
$(0:0:1)$. It is an amusing exercise to recover from this the $A_1$
singularity $cf=e^2$ at $P_x$ and the $A_2$ singularity $xdf=e^3$ at
$P_f$.

\section{Extrasymmetric format} \label{s!1}
\subsection{Extrasymmetric format}

Tom unprojections frequently lead to equations in extrasymmetric format.
Consider for example the $6\times6$ skew matrix\footnote{We omit the
diagonal terms (which are zero) and the $m_{ji} =-m_{ij}$ with $i < j$.}
\begin{equation}
N=
\begin{pmatrix}
z & y & a & b & d \\
& x & b & c & e \\
&& d & e & f \\
&&& \la z & \la y \\
&&&& \la x
\end{pmatrix} =
\begin{pmatrix} B & A \\ -A & \la B
\end{pmatrix}
\label{eq!3}
\end{equation}
A matrix of this shape is {\em extrasymmetric} (the term also covers
slightly more general cases, see \cite{TJ}, 9.1). It is made up of
$3\times3$ blocks, where the top right block $A$ is symmetric, the top
left block $B$ is skew, and the bottom right block $\la B$ repeats
the information contained in the top left block, in this case with a
scalar factor $\la$.

The $4\times4$ Pfaffians of $N$ generate the ideal of the Segre
embedding
\[
\Segre(\PP^2\times\PP^2)\subset\PP^8_{\Span{a,b,c,d,e,f,x,y,z}}.
\]
More precisely, the extrasymmetry means that the 15 upper-triangular
entries of $N$ consist of 9 independent entries and 6 repeats. The same
is true of the $4\times4$ Pfaffians of $N$, which give 9 relations and 6
repeats. The resulting 9 equations define a variety in
$\aff^9_{\Span{x,y,z,a,b,c,d,e,f}}$ that, for $\la\ne0$, is a linear
transformation away from the affine cone over
$\Segre(\PP^2\times\PP^2)$. The linear transformation involves taking
$\sqrt{-\la}$; swapping the signs of the square root interchanges
the two copies of $\PP^2\times\PP^2$. We leave the calculations as
entertainment.

A more banal way to define $\Segre(\PP^2\times\PP^2)$ is $\bigwedge^2
M=0$ with $M$ a generic $3\times3$ matrix. If we write $M=A+\sqrt{-\la}B$
with $A$ symmetric and $B$ skew, the ideal of $2\times2$ minors of $M$
equals the ideal of $4\times4$ Pfaffians of the extrasymmetric matrix
$N=\left(\begin{smallmatrix} B&A\\ -A&\la B\end{smallmatrix}\right)$.

More geometrically, this format displays $\PP^2\times\PP^2$ as a
nongeneric linear section of $\Grass(2,6)$.

\subsection{Specialise to $v_6(\PP(1,2,3))$}
Now we consider $\la$ as a variable and specialise the matrix
\eqref{eq!3} by setting $\la=0$, $z=0$ and $y=c$; the Pfaffian
equations specialise to \eqref{eq!2}. That is, the anticanonical cone
$X$ over $\PP(1,2,3)$ is the particular section $\la=0$, $z=0$ and $y=c$
of a degeneration of the cone over $\PP^2\times\PP^2$. Wiggling the
section gives one of the smoothing components of the deformations of
$X$.

\subsection{The same viewed as a Tom unprojection}

As we said, the extrasymmetric matrix $N$ in \eqref{eq!3} has 6 repeated
entries. The entries that are not repeated are the three diagonal
entries $a,c,f$ of the top right $3\times3$ block $A$. They correspond
to the three coordinate points of $\PP^2\times\PP^2$ such as
$P_a=(1:0:0;\ 1:0:0)$, etc. Here again $\la$ is a nonzero scalar.

Now project from $P_a$, and view the original equations as the result of
undoing this projection. A practical point of view on unprojection is
that it groups the 9 equations according to how they involve $a$.
Because of the format of \eqref{eq!3}, $a$ only appears linearly in 4
equations
\[
ac = \cdots,\quad ae = \cdots,\quad af = \cdots,\quad ax = \cdots,
\]
and the remaining 5 equations not involving $a$ are the Pfaffians of
\begin{equation}
\label{eq!4}
N_{\widehat4} =
\begin{pmatrix}
z& y& b& d \\
& x& c& e \\
&& e& f \\
&&& \la x
\end{pmatrix}
\end{equation}
(delete row and column 4 from $N$ of \eqref{eq!3}). What makes this a
Tom$_1$ matrix is that the 6 entries not in row and column 1 are in the
codimension 4 complete intersection ideal $(x,c,e,f)$. The coincidences
$m_{25}=m_{34}=e$ and $m_{45}=\la x=\la m_{23}$ that bring this about
are remnants of the extrasymmetry of $N$. From this point of view $a$ is
an unprojection variable, and the main theorem of \cite{PR} would allow
us to recover its equations.

Geometrically, the Pfaffians of \eqref{eq!4} define the projection of
$\PP^2\times\PP^2$ from $P_a$. It is a 4-fold section of $\Grass(2,5)$
containing the 3-plane $\PP^3_{\Span{b,d,z,y}}$ defined by the ideal
$(x,c,e,f)$.

\subsection{Finding the Tom format from $v_6(\PP(1,2,3))$}

We can start from the other end, dividing the 9 equations \eqref{eq!2} of
$v_6(\PP(1,2,3))$ into 4 that are linear in $a$ and 5 not involving
$a$. One gets $af=d^2$, $ae=bd$, $ac=b^2$ and $ax=bc$
together with the five Pfaffians of
\begin{equation}
\begin{pmatrix}
0& c& b& d \\
& x& c& e \\
&& e& f \\
&&& 0
\end{pmatrix}.
\label{eq!5}
\end{equation}
If we hope to describe the set of all 9 equations as
Pfaffians of a special $6\times6$ skew matrix, we must put $a$ where it
multi\-plies $x,c,e,f$ and {\em not} $b,d$, so put it as $m_{16}$.

\section{Double Jerry format} \label{s!2}

\subsection{Double Jerry} A neat starting point \cite[9.2]{TJ} is to
view Double Jerry as a theorem saying that a codimension 2 complete
intersection $m_1=m_2=0$ that contains two different codimension 3
complete intersections $(x_1,x_2,x_3)$ and $(y_1,y_2,y_3)$ is defined by
two bilinear forms
\begin{equation*}
m_1(x_1,x_2,x_3;\ y_1,y_2,y_3)
\quad\hbox{and}\quad m_2(x_1,x_2,x_3;\ y_1,y_2,y_3).
\end{equation*}
We can then introduce two parallel sets of unprojection equations
\begin{equation*}
s\cdot(x_1,x_2,x_3)=\cdots \quad\hbox{and}\quad t\cdot(y_1,y_2,y_3)=\cdots,
\end{equation*}
each taking us to codimension~3, together with a {\em long equation}
$st=\cdots$. Each unprojection separately is given by Cramer's rule,
leading to a $5\times5$ Pfaffian Jerry matrix, but the long equation is
an intriguing and in general surprisingly complicated function of
$m_1,m_2,x_i,y_i$. A particular case is worked out in Brown and
Georgiadis \cite{BG}.

\subsection{Our particular case} Rather than rework the general material
of \cite[9.2]{TJ}, consider only the case of the Newton polygon
\eqref{eq!1}.

As before, $a$ only appears linearly in 4 equations, so can be
eliminated or ``projected out'', expressing the variety as an
unprojection. The 5 equations not involving $a$ are again the Pfaffians
of \eqref{eq!5}. However, we now view it as a Jerry$_{23}$ matrix: in
fact, the 7 entries $\{0,x,c,e\}\cup \{c,x,e,f\}$ of its 2nd and 3rd
rows and columns consist of the regular sequence $x,c,e,f$ with repeats.
What makes it a {\em double} Jerry is that the pivot $m_{23}=x$ is one
of the variables on the nose, rather than a linear combination.

The matrix
\begin{equation}
\begin{pmatrix}
\mu f & c+\nu f & b & d  \\
& x & c & e \\
&& e & f \\
&&& -g
\end{pmatrix},
\label{eq!J45def}
\end{equation}
is a deformation respecting the Jerry$_{23}$ requirements just
described. Here $g$ is a new indeterminate of degree 1 and $\mu$ and
$\nu$ scalars. Putting back $a$ as unprojection variable defines a
family of del Pezzo 3-folds
\begin{equation*}
W_{\mu,\nu}\subset\PP^7_{\Span{a,b,c,d,e,f,g,x}}.
\end{equation*}
We recover $v_6(\PP(1,2,3))$ on setting $\mu=\nu=0$ and taking the
hyperplane section $g=0$.

\subsection{Interpretation as double Jerry}
\label{s!dJ}
Two of the Pfaffians of \eqref{eq!J45def} do not involve $x$:
\begin{equation}
be - cd + \mu fg \quad\hbox{and}\quad bf + cg - de + \nu fg
\label{eq!8}
\end{equation}
The codimension~2 complete intersection defined by these contains as
divisors two different codimension~3 complete intersection $V(b,d,g)$
and $V(c,e,f)$. Unprojecting these lead to $x$ and $a$ respectively.

In more detail, first write \eqref{eq!8} as
\[
\begin{array}r
\begin{pmatrix}
b & d & g
\end{pmatrix}
\\[9mm]
\end{array}
\kern-2mm
\begin{pmatrix}
e & f \\
-c & -e \\
\mu f & c+\nu f
\end{pmatrix}
=0.
\]
By Cramer's rule, $(b,d,g)$ is proportional to the minors of the
$3\times 2$ matrix. This predicts the remaining 3 minors of
\eqref{eq!J45def}:
\begin{equation}
\label{eq!11}
\begin{aligned}
xb &= c^2 - \mu ef + \nu cf, \\
xd &= ce + \nu ef - \mu f^2, \\
xg &= -cf + e^2.
\end{aligned}
\end{equation}

For $c,e,f$, working in the same way, \eqref{eq!8} gives
\[
\begin{array}r
\begin{pmatrix}
-d & b & \mu g \\
g & -d & b+\nu g
\end{pmatrix}
\\[7.5mm]
\end{array}
\kern-2mm
\begin{pmatrix} c \\ e  \\ f \end{pmatrix}
=0.
\]
Adjoining $a$ as the unprojection variable gives the other half of the
double Jerry:
\[
\begin{aligned}
ac &= b^2 + \nu bg + \mu dg, \\
ae &= bd + \nu dg + \mu g^2, \\
af &= -bg + d^2
\end{aligned} \quad\hbox{and}\quad
\begin{pmatrix}
-a & b & -d & \mu g \\
& -d & g & b+\nu g \\
&& f & -c \\
&&& e
\end{pmatrix}.
\label{eq!a}
\]
We get the long equation for $ax$ by cancelling $b,c,d,e,f$ or $g$ from
a linear combination of the other equations. There are many such
derivations: for example, start from $xg=e^2-fc$, multiply by $a$ and
rewrite the right hand side until it is divisible by $g$. The result is
\[
ax = (b+\nu g)(c+\nu f) - \mu(df - eg).
\]

The symmetry between the two unprojections is underlined by the fact
that the 9 equations are simply interchanged\footnote{They are also
invariant under $(d,e,\mu)\bij (-d,-e,-\mu)$. In these calculations
there may be several correct choices of signs (and many incorrect ones).
Getting the signs right can be a major headache, with no perfect
solutions.} by the involution
\[
\mu \bijj -\mu, \quad
a \bijj x, \quad b \bijj c, \quad d \bijj e, \quad f \bijj g.
\]

\subsection{$S_3$ symmetry} For general $\mu,\nu$, the 3-fold
$W_{\mu,\nu}$ is projectively equivalent to
$\PP^1\times\PP^1\times\PP^1$. Carrying this out requires an $S_3$
Galois field extension.

The little exercise in $A_2$ symmetry is fun and not quite obvious: the
three equations involving $x$ in \eqref{eq!J45def} are \eqref{eq!11}.
From them we deduce that in the deformation given by \eqref{eq!J45def},
the tag equation $xdf=e^3$ of \eqref{eq!1} deforms to
\begin{equation*}
x(df+eg)=e^3+\nu ef^2-\mu f^3=\Phi(e,f).
\end{equation*}

The projective equivalence of $W_{\mu,\nu}$ of
$\PP^1\times\PP^1\times\PP^1$ holds when the discriminant of $\Phi$ does
not vanish, and involves the roots of $\Phi$. It thus takes place over
its splitting field. The Galois group action permutes the 3 copies of
$\PP^1$. This reflects the Weyl group $W(A_2)=S_3$ symmetry behind the
deformation theory of the $A_2$ singularity.

Write $s,t,u$ for the roots of $\Phi$, so that $s+t+u=0$,
\begin{gather*}
\nu=st+ut+su=-(s^2+st+t^2), \quad \mu= ust=-st(s+t), \\
\hbox{and}\quad \Phi(e,f)=(e-sf)(e-tf)(e-uf),
\end{gather*}
Now set $y_0,y_1,y_2$ and $z_0,z_1,z_2$ to be the following
linear combinations of $(b,d,g)$ and $(c,e,f)$:
\[
\begin{aligned}
y_0 &= c+s\, e+tu\,f, \\
y_1 &= c+t\, e+su\,f,\\
y_2 &= c+u\,e+st\,f,
\end{aligned}
\quad\ \hbox{and}\quad\ 
\begin{aligned}
z_0 &=b-s \,d+tu\,g, \\
z_1 &=b-t\,d+su\,g, \\
z_2 &=b-u\,d+st\,g.
\end{aligned}
\]
After a calculation, we find
\[
\begin{aligned}
xz_i=y_jy_k,\\
ay_i=z_jz_k,\\
xa=y_iz_i,
\end{aligned}\qquad\hbox{for}\quad \{i,j,k\} = \{0,1,2\}.
\]
These are the standard
equations of $\PP^1\times\PP^1\times\PP^1$ as the $2\times2$ minors of
the 3-cube.

\section{Unprojection and deformations}
\label{s!3}

\subsection{Unprojection}
The general theory of unprojection was initiated by {\sc Kustin} and
{\sc Miller} \cite{KM} and developed in the present form by {\sc
Papadakis} and {\sc Reid}, see \cite{Ki, P, PR}.

Let $P \in D \subset X$ be a singular point of a Gorenstein scheme,
lying on a Gorenstein codimension~1 subscheme $D$. Consider the
adjunction sequence
\[
0 \to \om_X \to \Hom(\sI_D,\om_X) \to \om_D \to 0.
\]
By \cite[Lemma 1.1]{PR}, the $\Oh_X$-module $\Hom(\sI_D,\om_X)$ is
generated by two elements; we can take one of these as an injective
map $s\colon \sI_D\into \om_X\iso\Oh_X$ that projects to a basis element
of $\om_D\iso\Oh_D$. The {\em unprojection} $Y$ of $D$ in $X$ is the
spectrum of the $\Oh_X$-algebra $\Oh_X[S]/(S f_i-s(f_i))$, where the
$f_i$ generate the ideal $\sI_D\subset\Oh_X$. The scheme $Y$ is again
Gorenstein.

As $X$ is Gorenstein, $\Hom(\sI_D,\om_X)\iso \Hom(\sI_D,\Oh_X)$. We
calculate generators of $\Hom(\sI_D,\Oh_X)$ in concrete cases by
computer algebra, cf.~\cite{BP}. This construction also applies in a
relative situation, over a base space $T$. The most general $T$ is the
base of a versal deformation of the inclusion map $i\colon D \into X$.

\subsection{Combining the two deformation families}
In our case, the first order infinitesimal deformations of $i\colon D
\into X$ are described as the Pfaffian perturbations of the equations
contained in the ideal $(x,c,e,f)$. The trivial deformations are given
by vector fields $\Der(-\log D)$ preserving $D$. For deformations of
weight $-1$, this means that we make the matrix as general as possible,
with no coordinate transformations of $x$, $c$, $e$ and $f$ allowed. The
result is
\[
\begin{pmatrix}
z& c+y& b& d \\
& x& c& e \\
&& e& f \\
&&& -g
\end{pmatrix}
\]
The minus sign conforms with the deformation \eqref{eq!J45def}.

For deformations of weight $\ge0$ a short computation\footnote{available
from \url{http://www.math.chalmers.se/~stevens/singular.html}} in {\sc
Singular} \cite{DGPS} shows that the above deformations generate the
module of deformations: we can replace $y$ and $z$ with polynomials in
$f$, and $g$ with a polynomial in $x$ having deformation variables as
coefficients. Since our singularity is nonisolated some care is needed
with the meaning of infinite dimensional versal deformation. We restrict
ourselves here to deformations of nonpositive weight, that globalise to
deformations of the projective cone. Then the first order infinitesimal
deformations are given by
\begin{equation}
\begin{pmatrix}
z+\mu f& c+y+\nu f& b& d \\
& x& c& e \\
&& e& f \\
&&& -g+\la x
\end{pmatrix}
\label{eq!23def}
\end{equation}
For higher order deformations, the equations are the Pfaffians of the
matrix \eqref{eq!23def}, as the deformation is in particular a
deformation of $X$. The obstruction is that they must lie in the ideal
$(x,c,e,f)$. Hence setting these variables to zero in \eqref{eq!23def}
we find $gy=gz=0$ as the equations of the base space.

We compute $\Hom(\sI_D,\Oh_X)$ using {\sc Singular} \cite{DGPS} to
determine the unprojection, obtaining the equations
\begin{align*}
 & ac-b(b+\nu g)-\la(z+\mu f)^2-\mu dg+\la\nu c(c+y+\nu f), \\
 & ae-(b+\nu g)d-\la(c+y)(z+\mu f)-\mu g^2+\la \nu xd+\la\mu xg, \\
 & af-d^2+bg-\la(c+y)^2-\la\nu (c+y)f, \\
 & ax-(b+\nu g+\la\nu x)(c+y+\nu f)+d(z+\mu f)-\mu eg.
\end{align*}


We find two components, with total spaces isomorphic up to a smooth
factor with the Tom and Jerry formats of Sections \ref{s!1} and
\ref{s!dJ}. We replace $y$ by $y+c$ in the Tom equations, to obtain the
cone as section $\la=0$, $z=0$ and $y=0$. The coordinate transformations
needed are $a\mapsto a-\la\nu(c+y)$, $y\mapsto y-\nu f$ and $z\mapsto
z-\mu f$ for the Tom component and $a\mapsto a+\la\mu e$, $g\mapsto
g+\la x$ for Jerry. Note that these coordinate transformations mix the
deformation and the space variables.

\subsection{The versal deformation of the cone over $v_6(\PP(1,2,3))$}
Altmann \cite[Table 5.1]{A} records the result of our computation of
the infinite dimensional versal deformation.
What we have actually computed is the part in nonpositive weight, giving
the (embedded) versal deformation of the projective cone. After a simple
coordinate transformation and translation to our present coordinates,
the formulas there give exactly the same ideal as computed above in
terms of unprojection.

\subsection{The cone over an elliptic curve of degree 6}

The versal deformation of the cone over an elliptic normal curve of
degree~6 is described without equations by M\'erindol \cite{Me}. The
base space is the product of the cone over the Segre embedding of
$\PP^1\times\PP^2$ with the germ of an appropriate modular curve.

Deformations of negative weight can be described by Pinkham's
construction of ``sweeping out the cone''. More precisely, the total
space over a line in the base space is the cone over the anticanonical
model of an almost del Pezzo surface of degree 6, with the given
elliptic curve $E$ as hyperplane section. Such a surface is obtained by
blowing up three points on the curve, embedded in the plane by a linear
system of degree~3. M\'erindol's construction starts with a family of
such surfaces over an Abelian variety $A$, which is the hypersurface in
$\Pic^3\times E^3$ given by $3H-(P_1+P_2+P_3)=6O$. The Weyl group
$W=A_1\times A_2$ acts on this: $A_2$ permutes the three points, and
$A_1$ acts by
\begin{multline*}
(H;P_1,P_2,P_3)\mapsto \\
(2H-P_1-P_2-P_3; H-P_2-P_3; H-P_1-P_3,H-P_1-P_2).
\end{multline*}
Thus the base space of the versal deformation in negative weight is the
cone over $A/W\iso\PP^1\times \PP^2$.

We find the elliptic curve as hyperplane section of the singular del
Pezzo surface $v_6(\PP(1,2,3))$. In affine coordinates of $\PP^2$
related to \eqref{eq!dP} we take the curve $w^2=v^3+\ga v^2+v$,
realising the cone as the hyperplane section $f-x-\ga c-b=0$. Thus the
variable $a$ does not appear in the equation.

For the deformations of negative weight, we perturb the matrix
\eqref{eq!5} (with $b=f-x-\ga c$) with independent variables, subject to
the resulting equations lying in the ideal $(x,c,e,f)$. This means that
the entries multiplied by $m_{1,5}=d$ are not perturbed, and moreover,
no perturbation of $x$, $c$, $e$ or $f$ is absorbed by coordinate
transformations. We take
\begin{equation}
\begin{pmatrix}
z& c+y& b+u& d \\
& x& c& e+q \\
&& e& f+p \\
&&& s
\end{pmatrix}.
\label{eq!5def}
\end{equation}
The Pfaffians of this matrix with  $x$, $c$, $e$ and $f$ (and therefore
also $b$) equated to zero give the equations of the base space: the
minors of
\[
\begin{pmatrix}
z& y& u \\
q & p & s
\end{pmatrix}.
\]

As for the space of deformations of weight zero, a computation with {\sc Singular}
shows that it has dimension two. One deformation is given by the modulus
$\ga$, but there is another, corresponding to the choice of point from
which to project the curve.


The matrix \eqref{eq!5def} is neither a Tom nor a Jerry matrix. But it
can written in these forms after a small resolution of the base space.
We do this here for the Tom format. The cone over $\PP^1\times \PP^2$ is
resolved by $\PP^1\times \aff^3$. We introduce an inhomogenous
coordinate $\la$ on $\PP^1$ and set $q=\la z$, $s=\la y$ and $s=\la u$.
Then we can make the matrix into a Tom$_1$ by row and column operations.
After the coordinate transformation
\begin{multline*}
(c,d,e,f,x,u,y,z)\ \mapsto\ (c-\la x,d-\ga\la z+\la^2z,e, \\
f+\la c+\ga\la x-2\la^2x,x,u+b,y-c+\la x,z),
\end{multline*}
(so that $b=f-x-\ga c+\la c-2\la^2x$), the matrix takes the form
%
\[
\begin{pmatrix}
z& y& u& d \\
& x& c& e \\
&& e& f \\
&&& (\la+\ga\la^2+\la^3)x
\end{pmatrix}.
\]

\section{General remarks and FAQ}
\subsection{Which is Tom, and which is Jerry?}

We offer three answers as useful mnemonics. We do not assume any prior
familiarity with the Hanna--Barbera characters.
\begin{enumerate}
\renewcommand{\labelenumi}{(\roman{enumi})}
\item Tom is fatter. The ancestral Tom is the projective 4-fold
$\PP^2\times\PP^2$, whereas for Jerry it is the 3-fold
$\PP^1\times\PP^1\times\PP^1$.

\item The Tom$_i$ condition on a skew $5\times5$ matrix is that,
deleting the $i$th row and column, the remaining 6 entries $m_{jk}$ are
in a codimension~4 c.i.\ ideal. In simple cases, this means two
coincidences on the $m_{jk}$. On the other hand, the Jerry$_{jk}$
condition is that the 7 elements $m_{ij}=-m_{ji}$ and $m_{ik}=-m_{ki}$
in the $j$th and $k$th rows and columns are in a codimension 4 c.i.\
ideal, which means 3~conditions.

\item Weight-for-weight, Jerry is more singular. In fact any point
$P\in\PP^1\times\PP^1\times\PP^1\subset\PP^7$ lies on 3 lines, and the
linear projection from $P$ contracts these to nodes. In contrast, if we
take the 3-fold hyperplane section $V$ of $\PP^2\times\PP^2$ to get the
flag variety of $\PP^2$, the linear projection of $V$ from $P$ only has
two nodes.
\end{enumerate}

Trying to fit a Jerry unprojection into a $6\times6$ skew matrix format
is invariably a waste of time.

\subsection{What's it all about?}

A hypersurface or complete intersection is determined by the
coefficients of its defining equations, so its deformations are
unobstructed. The subtlety of the deformation theory in these cases is
nothing to do with obstructions, but how to pass to the quotient by the
appropriate equivalence relation, which involves dividing by the
groupoid of local diffeomorphisms.

The Buchsbaum--Eisenbud theorem \cite{BE} puts codimension~3 Gorenstein
ideals in the same framework: the variety is given by a skew
$(2k+1)\times(2k+1)$ matrix (most commonly $5\times5$), that encodes
both the defining equations and the syzygies, so that the entries of the
matrix can be freely deformed. In other words, the skew matrix is a
given mold, into which one can simply pour functions on the ambient
space in a liquid manner.

In contrast, one usually expects codimension~4 constructions to be
obstructed. A typical case is the cone over dP$_6$, whose deformation
theory has the 2 components we have mentioned many times.

The point of Tom and Jerry is that, in most commonly occurring cases, our
variety admits a Gorenstein projection to codimension~3, with the
projected variety given by the Pfaffians of a $5\times5$ skew matrix;
that is, the projected variety is a regular pullback from $\Grass(2,5)$
in its Pl\"ucker embedding, marked with an unprojection divisor that
corresponds to a linear subspace of $\Grass(2,5)$. Every geometer must
have done the easy exercise of seeing that any linear subspace of
$\Grass(2,4)$ (the Klein quadric) either consists of lines of $\PP^3$
passing through a point $P$, or dually, of lines contained in a plane
$\PP^2\subset\PP^3$. The Tom and Jerry formats answer the same question
for $\Grass(2,5)$; see \cite[2.1]{TJ}.

\subsection{Do they do everything?}
Unfortunately, no. Tom and Jerry provide two smooth components of the
deformation theory, and for deformation problems entirely contained
within one component or the other, they can be relied on to do
everything. However, we know other cases in codimension~4 that appear
not to have any useable structure of Kustin--Miller unprojection.

A general structure theorem for Gorenstein codimension 4 ideals is
described in \cite{G4}. It is rather complicated, as it should account
for the singular total spaces of versal deformations, cf.~the discussion
of the cone over an elliptic curve of degree 6 above. Deformations of
its hyperplane sections are even more complicated.


\end{document}